\input amstex
\documentstyle{amsppt}
\pageheight{204mm}
\pagewidth{133mm}
\magnification\magstep1


\def\nologo{\let\logo@\empty}

\def\Aut{\operatorname{Aut}}
\def\Aut{\operatorname{Aut}}

\def\End{\operatorname{End}}

\def\gr{\operatorname{gr}}

\def\Hom{\operatorname{Hom}}

\def\LMH{\operatorname{LMH}}

\def\rank{\operatorname{rank}}

\def\SL{\operatorname{SL}}

\def\toric{\operatorname{toric}}
\def\abtoric{\operatorname{|toric|}}
\def\torus{\operatorname{torus}}
\def\abtorus{\operatorname{|torus|}}

\def\tCu{\tsize\bigcup}

\def\bC{\bold C}

\def\bN{\bold N}

\def\bQ{\bold Q}

\def\bR{\bold R}

\def\bZ{\bold Z}

\def\cB{{\Cal B}}
\def\cC{{\Cal C}}

\def\cO{{\Cal O}}

\def\fg{{\frak g}}

\def\G{\Gamma}

\def\sig{\sigma}
\def\Sig{\Sigma}

\def\.{$.\;$}

\def\gp{{\text{\rm gp}}}

\def\mult{{\text{\rm mult}}}

\def\resp.{\text{\rm resp}.\;}

\def\val{{\text{\rm val}}}

\def\O^logten{\cO\^log\otimes}

\let\bs=\backslash

\let\lan=\langle
\let\lan=\langle

\let\ran=\rangle

\let\sub=\subset

\def\Dc{\check{D}}
\def\Ec{\check{E}}

\def\LMH{\text{LMH}}

\topmatter

\title
Analyticity of the closures of some Hodge theoretic subspaces
\endtitle

\author
Kazuya Kato\footnote{\text{Partially supported by NFS grant DMS 1001729.}},
Chikara Nakayama\footnote{\text{Partially supported by JSPS Grants-in-Aid for Scientific Research (C) 18540017,
(C) 22540011.}},
Sampei Usui\footnote{\text{Partially supported by JSPS Grant-in-Aid for Scientific Research (B) 19340008.}}
\endauthor

\address
\newline
{\rm Kazuya KATO}
\newline
Department of Mathematics
\newline
University of Chicago
\newline
5734 S.\ University Avenue
\newline
Chicago, Illinois 60637, USA
\newline
{\tt kkato\@math.uchicago.edu}
\endaddress

\address
\newline
{\rm Chikara NAKAYAMA}
\newline
Graduate School of Science and Engineering
\newline
Tokyo Institute of Technology 
\newline
Meguro-ku, Tokyo, 152-8551, Japan
\newline
{\tt cnakayam\@math.titech.ac.jp}
\endaddress

\address
\newline
{\rm Sampei USUI}
\newline
Graduate School of Science
\newline
Osaka University
\newline
Toyonaka, Osaka, 560-0043, Japan
\newline
{\tt usui\@math.sci.osaka-u.ac.jp}
\endaddress

\abstract
In this paper, we prove a general theorem concerning the analyticity of 
the closure of a subspace defined by a family of variations of mixed 
Hodge structures, which includes 
the analyticity of the zero loci of degenerating normal functions.
  For the proof, we use a 
moduli of the valuative version of log mixed Hodge structures. 
\endabstract

\footnote"{}"{2010 {\it Mathematics Subject Classification}.
Primary 14C30; Secondary 14D07, 32G20.}

\NoRunningHeads

\endtopmatter

\document

\head
\S0. Introduction
\endhead

 We recall the definition of the classifying space $D$ 
of mixed Hodge structures 
with polarized graded quotients, introduced in \cite{U84}, which 
is a natural generalization of Griffiths classifying space of 
polarized Hodge structures (\cite{G68}).

  Fix a $4$-ple $\Lambda:=
(H_0, W, (\lan\;,\;\ran_k)_k$, $(h^{p,q})_{p,q}),
$
where 
$H_0$ is a finitely generated free $\bZ$-module, 
$W$ is an increasing filtration on $H_{0,\bQ}:=\bQ\otimes_\bZ H_0$, 
$\lan\;,\;\ran_k$ is a non-degenerate $\bQ$-bilinear form $\gr^W_k \times$ $
\gr^W_k \to \bQ$ given for each $k\in \bZ$ which is symmetric if $k$ is even and anti-symmetric if $k$ is odd, and 
$h^{p,q}$ is a non-negative integer given for $p, q\in\bZ$ such
that $h^{p,q} = h^{q,p}$ for all $p, q$, and that 
$
\rank_\bZ(H_0) = \tsize\sum_{p,q} h^{p,q}$, $ 
\dim_\bQ(\gr^W_k)= \tsize\sum_{p+q=k} h^{p,q} \;\;\text{for all} \;k.
$

For $A=\bZ, \bQ, \bR,$ or $\bC$, let
$G_A$ be the group of all $A$-automorphisms of $H_{0,A}$ which are 
compatible with 
$W$ and $\lan\;,\;\ran_k$ for any $k$. 
For $A=\bQ, \bR$, or $\bC$, let
$\fg_A$ be the Lie algebra of $G_A$, which is identified with a 
subalgebra of $\End_A(H_{0,A})$ in the natural way 
(cf.\ \cite{KNU10} 1.2).

Let $D$ be the set of all decreasing filtration $F$ on $H_{0,\bC}$ for which $(H_0, W, F)$ is a mixed Hodge structure such that the $(p, k-p)$ Hodge number of $F(\gr^W_k)$
coincides with $h^{p,k-p}$ for any $p, k\in \bZ$ and such that $F(\gr^W_{k})$ is polarized by $\lan\;,\;\ran_k$ for all $k$. 
  Let $\check D$ be its compact dual (cf.\ \cite{KNU10} 1.3). 

In this paper, we construct valuatively toroidal partial compactifications of 
$D$.  As an application, we prove the following theorem.

\proclaim{Theorem 0.1} 
Let $S$ be a complex analytic space, let $Y$ be a closed analytic subspace of $S$, and let $S^*=S-Y$. 
Assume that $S^*$ is non-singular and is dense in $S$.  
Let $\G$ be a neat subgroup of $G_\bZ$. For $1\leq j\leq n$, let $f_j: S^*\to \G \bs D$ be the period map associated to a variation of mixed Hodge structure with polarized graded quotients which is admissible with respect to $S$ {\rm(\cite{K86}).} 
Let $Z^*=\{s\in S^*\;|\; f_1(s)=\dots=f_n(s)\}$, and let $\overline {Z^*}$ be the closure of $Z^*$ in $S$. Then $\overline{Z^*}$ is an analytic subset of $S$.  
\endproclaim 

  We can deduce from Theorem 0.1 the following theorem 
proved by Brosnan--Pearlstein and, independently, by Schnell.

\proclaim{Theorem 0.2 {\rm(Brosnan--Pearlstein [BP.p], Schnell [Scn.p])}} 
  Let $S^*$ be a smooth complex algebraic variety. 
  Let $H$ be a variation of PHS of weight $-1$ on $S^*$. 
  Let $\nu:S^* \to J(H)$ be an admissible normal function {\rm (\cite{Sa96})}, 
where $J(H)$ is the intermediate Jacobian. 
  Then the zero locus $Z(\nu) \sub S^*$ of $\nu$ is algebraic. 
\endproclaim

\definition{Remark}
  The case of 0.1 where $S$ is smooth and $Y$ is a smooth divisor 
is also proved in \S4 in \cite{KNU10}. 
  The case of 0.2 where $S^*$ is 
compactifiable by a smooth divisor 
was also proved by 
Brosnan--Pearlstein (\cite{BP09a}, \cite{BP09b}) 
and, independently, by Saito (\cite{Sa.p}).
\enddefinition

  To deduce 0.2, after taking a smooth compactification of $S^*$, 
we apply the case of 0.1 where $\gr^W_k=0$ unless $k=0, 1$ and 
$\gr^W_0=\bQ$.  
  See Remark 1 of \S4 in \cite{KNU10}.

  After a preparation in \S1, we state some theorems on moduli spaces 
in \S2--\S3. 
  We prove them and some more propositions in \S4, 
which are essentially reduced to our previous works. 
  We apply them to prove Theorem 0.1 in \S5--\S6.
  
\bigskip

\head
\S1. Variants of toric varieties
\endhead

We recall some facts about toric varieties in the form which will be used
 in this paper.

\medskip

{\bf 1.1.}  
Fix a subgroup $\Gamma$ of $G_\bZ$. Let $\sigma$ be a sharp rational nilpotent cone (\cite{KNU10} 2.1.1) 
in $\fg_\bR$. Let  $\G(\sigma):= \G\cap \exp(\sigma)$
and assume that $\sig$ is generated by $\log(\G(\sig))$ as a cone. Let 
 $P(\sigma):= \Hom(\G(\sigma), \bN)$, the dual monoid.
Define the toric variety and the torus associated to $\sigma$ by
$$
\toric_\sigma:= \Hom(P(\sigma), \bC^{\mult}) \supset 
\torus_\sigma:= \Hom(P(\sigma), \bC^\times)=\G(\sig)^{\gp}\otimes_{\bZ} \bC^\times, 
$$
where $\bC^{\mult}$ denotes $\bC$ regarded as a multiplicative monoid.

\medskip

{\bf 1.2.}  
We continue to use the previous notation.

Let $\toric_{\sig,\val}:=(\toric_{\sig})_{\val}$ be the projective limit of the log 
modifications of $\toric_{\sig}$, that is, the projective limit of toric varieties over $\bC$ associated to rational finite subdivisions of the cone $\sig$.

\medskip

{\bf 1.3.}  
$\abtoric_\sigma$ and $\abtorus_\sigma$ are defined as $\Hom(P(\sigma), \bR_{\ge0}^{\mult})$ and $\Hom(P(\sigma), \bR_{>0}^{\mult})=\Gamma(\sig)^{\gp}\otimes_{\bZ} \bR_{>0}^{\mult}$, respectively (\cite{KNU10} 2.2.3).

$\abtoric_{\sig,\val}$ is defined as the closure of the inverse image of $\abtorus_\sigma$ in $\toric_{\sig,\val}$.

\medskip

{\bf 1.4.}  
Summing up the above, we have a commutative diagram
$$
\matrix
&\torus_{\sig}&\hookrightarrow & \toric_{\sig,\val}&\hookleftarrow &\abtoric_{\sig,\val} \\
&\shortparallel &&\downarrow &&\downarrow \\
&\torus_{\sig}&\hookrightarrow & \toric_{\sig\phantom{,\val}}&\hookleftarrow&\abtoric_{\sig}.{}_{\phantom{\val}} \\
\endmatrix
$$

\bigskip

\head
\S2. $D_{\val}$ and $ D_{\val}^{\sharp}$
\endhead

  In this section, we expose some general theory on the valuatively toroidal partial compactifications of $D$. 
  We will develop this theory in detail 
in forthcoming parts of our series of papers 
on classifying spaces of degenerating mixed Hodge structures 
(\cite{KNU09}, \cite{KNU11}, \cite{KNU.p}, ...). 
  Here we state and prove the part which suffices for the proof of Theorem 0.1, but only outline the remaining part of the whole theory. 
  
\medskip

{\bf 2.1.}  
We recall the definition of the sets $D_{\val}$ and $ D_{\val}^{\sharp}$ in \cite{KNU.p}  4.2.1.

Let $D_{\val}$ (resp. $D_{\val}^{\sharp}$) be the set of all triples $(A, V, Z)$, where
$A$ is a $\bQ$-subspace of $\fg_\bQ$ consisting of mutually commutative nilpotent elements, 
$V$ is a sharp valuative submonoid of $A^*:=\Hom_\bQ(A, \bQ)$, and 
$Z$ is a subset of $\Dc$ such that $Z=\exp(A_\bC)F$ (resp.\ $Z=\exp(iA_\bR)F$) for any $F\in Z$ and that there exists a finitely generated rational subcone $\tau$ of $A_\bR$ satisfying the conditions that
$Z$ is a $\tau$-nilpotent orbit (resp. $i$-orbit) and $(A \cap \tau)^\vee \subset V$ in $A^*$. Here $A_\bR:= \bR\otimes_{\bQ} A\subset \fg_\bR$ and $A_\bC:=\bC\otimes_{\bQ} \fg_\bC$.

We have a canonical surjective map $D_{\val}^{\sharp}\to D_{\val}$, $(A,V,Z) \mapsto (A,V, \exp(A_\bC)Z)$. 
\medskip

{\bf 2.2.}  
We endow $D_{\val}^\sharp$ with the strongest topology for which the map $D_{\sig,\val}^{\sharp}\to D_{\val}^{\sharp}$ is continuous for any sharp rational nilpotent cone 
$\sig$. 
  Here the topology of $D_{\sig,\val}^{\sharp}$ is defined in \cite{KNU.p} 
4.2.7. 

\proclaim{Theorem 2.3} {\rm (i)} $D_{\val}^{\sharp}$ is Hausdorff. 

\medskip

{\rm (ii)} Let $\G$ be a subgroup of $G_{\bZ}$. 
Then the action of $\G$ on $D_{\val}^{\sharp}$ is proper. The quotient $\G \bs D_{\val}^{\sharp}$ is Hausdorff. 
If $\G$ is neat, the projection $D_{\val}^{\sharp} \to \G \bs D_{\val}^{\sharp}$ is a local homeomorphism.
\endproclaim

\medskip

{\bf 2.4.} 
Let $\G$ be a subgroup of $G_\bZ$. 
Let $\cC_{\G}$ be the set of all 
sharp rational nilpotent cones in $\fg_\bR$ generated by the logarithms of a finite number of elements of $\G$. 
  Let $D_{\val,(\G)} \sub D_{\val}$ be the union of $D_{\sig,\val}$ for 
$\sig \in \cC_{\G}$. 

We endow $\G \bs D_{\val,(\G)}$ with 
the topology, the sheaf of (complex) analytic functions and  the log structure 
as follows. 

A subset $U$ of $\G \bs D_{\val,(\G)}$ is open if and only if, for any 
$\sig\in \cC_{\G}$, the inverse image of $U$ in $\G(\sig)^{\gp} \bs D_{\sig,\val}$ is open. 
For an open set $U$ of $\G \bs D_{\val,(\G)}$, a complex valued function $f$ on $U$ is analytic if and only if, for any $\sig \in \cC_{\G}$, the pull-back of $f$ on the inverse image of $U$ in $\G(\sig)^{\gp}\bs D_{\sig,\val}$ is analytic.
   The log structure of $\G\bs D_{\val,(\G)}$ is the following subsheaf of the sheaf of analytic functions on $\G \bs D_{\val,(\G)}$. For an open set $U$ of $\G \bs D_{\val,(\G)}$ and for an analytic function $f$ on $U$, $f$ belongs to the log structure of $\G \bs D_{\val,(\G)}$ if and only if, for any $\sig\in\cC_{\G}$, the 
pull-back of $f$ on the inverse image of $U$ in 
$\G(\sig)^{\gp}\bs D_{\sig,\val}$ belongs to the log structure of $\G(\sig)^{\gp}\bs D_{\sig,\val}$. 

  In the above, the structure of a log local ringed space 
$\G(\sig)^{\gp}\bs D_{\sig,\val}$ is defined in \cite{KNU.p} 4.2.7. 

  In a forthcoming paper, we will prove nice properties of 
$\G \bs D_{\val,(\G)}$ which are analogous to Theorem 2.3. 

\medskip

{\bf 2.5.}
To state the next theorem, 
let $\G$ be a neat subgroup of $G_{\bZ}$.
We denote by $\LMH_{(\Lambda,\G),\val}$ the sheaf on 
the category $\cB(\log)$ (\cite{KU09} 3.2.4) with the usual topology 
associated to the presheaf 
sending $S \in \cB(\log)$ to 
the inductive limit of the sets of 
all isomorphism classes of log mixed Hodge structures on $S'$ with polarized graded quotients of the given Hodge type $\Lambda$ 
and with $\G$-level structure, where $S'$ runs over the sets of log 
modifications of $S$ (\cite{KU09} 3.6). 

  Let $S$ be an object of $\cB(\log)$. 
  A morphism $f: S_{\val} \to \G \bs D_{\val,(\G)}$ of log local ringed spaces over $\bC$ is 
called a {\it good morphism}  if the following condition (1) is satisfied.

\medskip

  (1) For any point $p$ of $S_{\val}$, there are an open neighborhood $U$ of $p$ in $S_\val$, $\sig \in \cC_{\G}$, and a morphism $U \to \G(\sig)^{\gp} \bs D_{\sig,\val}$ which induces 
the restriction $f\vert_U$ such that the composite $U \to 
\G(\sig)^{\gp} \bs D_{\sig,\val} \to 
\G(\sig)^{\gp} \bs D_{\sig}$ factors through an open subspace $V$ of a 
log modification of an open subspace of $S$. 
$$
\matrix
S_\val&\supset&U&@>>>&V\\
\\
@V{f}VV@VVV@VVV\\
\\
\G\bs D_{\val,(\G)}&@<<<&\G(\sigma)^\gp\bs D_{\sigma,\val}&@>>>&\;\G(\sigma)^\gp\bs D_\sigma.
\endmatrix
$$

\proclaim{Theorem 2.6} 
Let $\G$ and $S$ be as above.
Then there is a canonical functorial bijection between 
$\text{\rm{LMH}}_{(\Lambda,\G),\val}(S)$
and 
the set of good morphisms $S_{\val} \to \G \bs D_{\val,(\G)}$.
\endproclaim

  We expect that 
all morphisms $f\colon S_{\val} \to \G \bs D_{\val,(\G)}$ are good 
(hence $\G \bs D_{\val,(\G)}$ is in fact a fine moduli), though we have not yet proved it.

  Notice that in analogous 
theorems on moduli of (mixed) log Hodge structures in \cite{KU09} and \cite{KNU.p}, we always fix a fan or a weak fan $\Sig$ to restrict local monodromies.
Thus, in those works, there is always a problem of different nature 
to make such a 
$\Sig$.
Here in 2.6, we do not need any $\Sig$ in virtue of the valuative 
formulation.  

In the proof of Theorem 0.1 below, we only use a part of Theorem 2.6, but we expect that 2.6 
should relate more closely to 0.1 and that 2.6 will play more important roles 
in some future applications.

The proofs of these theorems 2.3 and 2.6 will be given in \S4.

\bigskip

\head
\S3. $D^{\sharp}_{(\sig),\val}$
\endhead

{\bf 3.1.} 
To prove Theorem 2.3 and Theorem 2.6, we introduce the following set
$D_{(\sig),\val}^{\sharp}$, where $\sig$ is a sharp rational nilpotent cone. 

Let 
$$
D^{\sharp}_{(\sig),\val}=\tsize\bigcup_{\tau} D^{\sharp}_{\tau,\val}\subset 
D^{\sharp}_{\val},
$$
where $\tau$ ranges over all rational subcones of $\sig$. 

\medskip

{\bf 3.2.} 
We endow $D_{(\sig),\val}^\sharp$ with the strongest topology for which the map $D_{\tau,\val}^{\sharp}\to D_{(\sig),\val}^{\sharp}$ is continuous for any rational subcone 
$\tau$ of $\sig$. 

\medskip

{\bf 3.3.} 
For a sharp rational nilpotent cone $\sig$, we define the 
topological space 
$E_{(\sig),\val}^{\sharp}$ as follows. 

Let $\Ec^{\sharp}_{\sig,\val} = \abtoric_{\sigma,\val} \times \Dc$.
  Recall that we defined $E^{\sharp}_{\sig, \val}:= \abtoric_{\sigma,\val} \times_{\abtoric_\sigma} E^{\sharp}_{\sig} \sub \Ec^{\sharp}_{\sig,\val}$ in \cite{KNU.p} 4.2.5.

  We define the set 
$E^{\sharp}_{(\sig),\val}$ such that 
$E^{\sharp}_{\sig, \val}\subset E^{\sharp}_{(\sig),\val}\sub \Ec^{\sharp} _{\sig,\val}$ as follows. 
  Let $\tau$ be a rational subcone of $\sig$. 
  Let $U(\tau)$ be the space 
$\abtorus_{\sig} \times^{\abtorus_{\tau}}E^{\sharp}_{\tau,\val}$. 
  Then we have the map $U(\tau) \to \Ec^{\sharp} _{\sig,\val}$ 
induced by the map 
$\abtorus_{\sig} \times^{\abtorus_{\tau}}\abtoric_{\tau,\val} \to 
\abtoric_{\sig,\val}$.
  We define $E^{\sharp}_{(\sig),\val}$ as the union of the images of all such 
$U(\tau) \to \Ec^{\sharp} _{\sig,\val}$, where $\tau$ runs over the set of the rational 
subcones of $\sig$. 

  We endow $E^{\sharp}_{(\sig),\val}$ with the strongest topology for which 
the maps $U(\tau) \to E^{\sharp}_{(\sig),\val}$ is continuous for any 
rational subcone $\tau$ of $\sig$.

We have a canonical projection
$
E_{(\sig),\val}^{\sharp}\to D_{(\sig),\val}^{\sharp}.
$

Let $\sig_\bR\sub\fg_\bR$ be the $\bR$-linear span of $\sig$. The continuous action of $i\sig_\bR$ on $E_{\sig,\val}^{\sharp}$ in \cite{KNU.p} 5.2.1 is naturally extended to a continuous action of $i\sig_\bR$ on $E_{(\sig),\val}^{\sharp}$.
\medskip

\proclaim{Theorem 3.4} {\rm (i)}
The action of $i\sig_{\bR}$ on $E_{(\sig),\val}^{\sharp}$ is free and proper. 
\medskip
{\rm (ii)} $E_{(\sig),\val}^{\sharp}$ is an $i\sig_\bR$-torsor over 
$D_{(\sig),\val}^{\sharp}$ in the category of 
topological spaces.
\endproclaim

\proclaim{Theorem 3.5} 
The inclusion $D_{(\sig),\val}^{\sharp} \to D_{\val}^{\sharp}$ is an open map.
\endproclaim 

The proofs of these theorems will be given in \S4.

\bigskip

\head
\S4. Proofs of the theorems in \S2 and \S3
\endhead

The proofs of the theorems in \S2--\S3 go on in a similar way as those of the  
theorems in \cite{KNU.p}.

\medskip

{\bf 4.1.} {\it Proof of Theorem 3.4.}
  The freeness in (i) is reduced to \cite{KNU.p} 5.2.2 (ii). 
  The properness in (i) is an analogue of ibid.\ 5.2.5, and 
is similarly reduced to the variant of ibid.\ 5.2.3 (ii) with 
the subscripts \lq\lq$(\sig)$'' instead \lq\lq$\sig$''.
  The proof of this variant is the same. 
Note that, by the definition of the topology of $D_\val^\sharp$ in 2.2, the map $\psi : D_\val^\sharp \to D_{\SL(2)}^I$ in \cite{KNU.p} 4.3 is continuous.
Note also that the topology of $E^{\sharp}_{(\sig),\val}$ given in 3.3 
is stronger or equal to the strong 
topology as a subset of $\Ec^{\sharp}_{\sig,\val}$.

  The assertion (ii) is an analogue of \cite{KNU.p} 5.2.8. 
  By ibid.\ 5.2.7 (a) and by the assertion (i), it is enough to 
check the condition (1) in ibid.\ 5.2.7
(with $H=i\sig_{\bR}$ and $X=E^{\sharp}_{(\sig),\val}$), 
which is done similarly as follows. 
  Let $x=(q,F) \in E^{\sharp}_{(\sig),\val}$. 
  Let $A \sub \Ec_{\sig}$ be the analytic space containing the image of $x$ 
constructed in the same way as in the pure case \cite{KU09} 7.3.5. 
  Let $S_1$ be the pull-back of $A$ in $E^{\sharp}_{(\sig),\val}$, 
let $U$ be a sufficiently small neighborhood of $0$ in $\sig_{\bR}$, 
and let $S=\{(q',\exp(a)F')\,|\,(q',F') \in S_1, a \in U\}$. 
  Then, $S$ satisfies the condition (1) in ibid.\ 5.2.7. 

\medskip

{\bf 4.2.} {\it Proof of Theorem 3.5.}
This is an analogue of \cite{KNU.p} 5.3.1. 
In the same way as there, we reduce 3.5 to the following (1). 

(1) If $\sig$ is a sharp rational nilpotent cone and $\tau$ is a 
rational subcone, then the inclusion map $D^{\sharp}_{(\tau),\val} \to 
D^{\sharp}_{(\sig),\val}$ is an open map. 

The proof of (1) is also similar.  We use the fact that the topology 
of $D^{\sharp}_{(\sig),\val}$ is the quotient topology from 
$E^{\sharp}_{(\sig),\val}$. 

\medskip

{\bf 4.3.} {\it Proof of Theorem 2.3.}
  The assertion (i) is an analogue of \cite{KNU.p} 5.3.3.
  By 3.5, it is sufficient to prove the variant of (1) 
in the proof of ibid.\ 5.3.3  with 
the subscripts \lq\lq$(\sig)$'' by \lq\lq$\sig$''.
  Similarly to ibid.\ 5.3.3, by 3.4 (ii) (which is proved in 4.1), this 
variant of ibid.\ 5.3.3 (1) is reduced 
to the variant of ibid.\ 5.2.3 (i) used in 4.1.

  The assertion (ii) is an analogue of \cite{KNU.p} 5.3.6. 
As we have seen in 4.1, the map $\psi : D_\val^\sharp \to D_{\SL(2)}^I$ is continuous.
By \cite{KNU11}  3.5.17, the action of $\G$ on $D_{\SL(2)}^I$ is proper. 
By \cite{KNU.p} 5.2.4.3 (ii), the properness of the action of $\G$ on 
$D_\val^\sharp$ follows from these two facts together with 2.3 (i).
  Then, $\G \bs D^{\sharp}_{\val}$ is Hausdorff by ibid.\ 5.2.4.1.
By ibid.\ 5.2.4.2, the last assertion follows from the above result and the freeness of the action of $\G$ on $D_\val^\sharp$ which is easily derived from ibid.\ 5.3.5 (i).

\medskip 

The next is a complement to \cite{KNU.p}. 

\proclaim{Proposition 4.4} 
Let $\G$ be a subgroup of $G_{\bZ}$ and 
let $\Sig$ be a fan in $\fg_{\bQ}$ {\rm (\cite{KNU10} 2.1.2)} 
which is strongly compatible with $\G$ {\rm (\cite{KNU10} 2.1.5)}. 
Then the following holds. 

\medskip

{\rm (i)}
There is a canonical isomorphism 
$(\G \bs D_{\Sig})_{\val} = \G \bs D_{\Sig,\val}$ 
of log local ringed spaces over $\bC$. 

\medskip

{\rm (ii)} 
There is a canonical homeomorphism 
$(\G \bs D_{\Sig})_{\val}^{\log} = \G \bs D_{\Sig,\val}^{\sharp}$ 
of topological spaces. 
\medskip

For the definitions of the right hand sides, see \cite{KNU.p} {\rm 4.2.7.}
\endproclaim

\demo{Proof}
  Let $\sig \in \Sig$, and we may replace $\G$ by $\G(\sig)^{\gp}$ and 
$\Sig$ by $\sig$. 

We prove (i).  Let $\sig_\bC\sub \fg_\bC$ be the $\bC$-linear span of $\sig$. By \cite{KNU.p} 2.5.3, $\G(\sig)^{\gp} \bs D_{\sig}$ is 
the quotient $\sig_{\bC} \bs E_{\sig}$. 
  By this, we see that $(\G(\sig)^{\gp} \bs D_{\sig})_{\val}$ is 
the quotient $\sig_{\bC} \bs E_{\sig,\val}$. 
  On the other hand, it is similar to ibid.\ 2.5.3 to see that 
$\G(\sig)^{\gp} \bs D_{\sig,\val}$ is also the quotient 
$\sig_{\bC} \bs E_{\sig,\val}$ (cf.\ ibid.\ 5.4.6). 
  Hence we have (i). 

  The assertion (ii) is similarly proved by using \cite{KNU.p} 5.2.8 instead 
of ibid.\ 2.5.3. 
\qed
\enddemo

{\bf 4.5.} {\it Proof of Theorem 2.6.}

This is an analogue of \cite{KNU.p} 2.6.6 (whose 
pure Hodge theoretic version is \cite{KU09} Theorem B), and is reduced to 
\cite{KNU.p} 2.6.6 as follows.

Assume that we are given an element $H$ of 
$\LMH_{(\Lambda,\G),\val}(S)$. 
  We provide the corresponding good morphism $S_{\val} \to \G \bs D_{\val,(\G)}$. 
  Assume first that $H$ comes from an LMH over $S$.
  Then, by the mixed Hodge theoretic version of \cite{KU09} 4.3.8 (which is used in the proof of \cite{KNU.p} 7.4.1) and by 
\cite{KNU.p} 2.6.6, locally on $S$, there 
are a log modification $S'$, a subgroup $\G'$ of $\G$, and a fan $\Sig$ 
in $\fg_{\bQ}$ which is strongly compatible with $\G'$ such that 
$H$ comes from a morphism $S' \to \G' \bs D_{\Sig}$. 
  Taking $(-)_{\val}$, we have $(S')_{\val} \to (\G' \bs D_{\Sig})_{\val}=
\G' \bs D_{\Sig,\val} \to \G \bs D_{\val,(\G)}$ by 4.4 (i).
In the general case, the above construction glues to give the desired morphism. 

  Conversely, let $S_{\val} \to \G \bs D_{\val,(\G)}$ be a good morphism.
  Then, again by \cite{KNU.p} 2.6.6, locally on $S$, we have LMHs of type 
$\Lambda$ with $\G$-level 
structures on various open subspaces of various log modifications.
  Since $S_{\val} \to S$ is proper, these LMHs glue into an LMH on 
one log modification locally on $S$. 
Thus we get an element of $\LMH_{(\Lambda,\G),\val}(S)$. 
\qed

\medskip

  The next is proved by the same argument as in 4.5, by using 4.4 (ii) instead of 4.4 (i), and will be used later in 6.1 for the proof of Theorem 0.1.

\proclaim{Proposition 4.6} 
Let $S$ be an object in $\cB(\log)$. 
Let $\G$ be a neat subgroup of $G_{\bZ}$.
Let $H \in \text{\rm{LMH}}_{(\Lambda,\G),\val}(S)$.
Then there is a natural continuous map 
$S_{\val}^{\log} \to \G \bs D_{\val}^\sharp$ which is 
compatible with the good morphism 
$S_{\val} \to \G \bs D_{\val,(\G)}$ corresponding 
to $H$ by $2.6$.
\endproclaim 

\demo{Proof}
The proof goes exactly in parallel to the former half of 4.5.
Let $S' \to \G' \bs D_{\Sig}$ be the morphism in 4.5 from which $H$ comes locally.
Taking $(-)_{\val}^{\log}$ this time, we have $(S')_{\val}^{\log} \to (\G' \bs D_{\Sig})_{\val}^{\log} = \G' \bs D_{\Sig,\val}^\sharp \to \G \bs D_{\val}^\sharp$ by 4.4 (ii).
  The rest is the same. 
\qed
\enddemo

\bigskip

\head
\S5. Preparation for the proof of Theorem 0.1
\endhead

{\bf 5.1.} 
Let the situation be as in the hypothesis of Theorem 0.1.

By resolution of singularity, we may and do assume that $S$ is non-singular and $Y$ is a divisor with normal crossings. 
We endow $S$ with the log structure associated to $Y$.

To prove Theorem 0.1, it is sufficient to show the following assertion (A).

\medskip

(A) Let $s\in S$. 
Then there is an open neighborhood $U$ of $s$ in $S$ such that 
$\overline{Z^*} \cap U$ is an analytic subset of $U$.

\medskip

We replace $S$ by a small open neighborhood of $s$ in $S$.
Replacing $S$ further by a finite ramified covering, we may and do assume that the local monodromy groups along $Y$ are unipotent.

\medskip

{\bf 5.2.}
Let $H^{(j)}$ be the variation of mixed Hodge structure on $S^*$ corresponding to $f_j$. 
  Then, $H^{(j)}$ extends uniquely to a log mixed Hodge structure on $S$, which we still denote by $H^{(j)}$. 

  We explain this extension. 
Let $\tau: S^{\log}\to S$ be as in \cite{KU09} 0.2.9.
  First, the lattice $H^{(j)}_{\bZ}$ extends as a locally 
constant sheaf over $S^{\log}$, which we still denote by the same symbol 
and which is isomorphic modulo $\G$ to the constant sheaf $H_0$.
  Further, by the admissibility and 
by the nilpotent orbit theorem of Schmid (\cite{Scm73}) reformulated 
as in \cite{KU09} 2.5.13--2.5.14 on each graded quotient $\gr^W_w$, $H^{(j)}$ extends uniquely to a log mixed Hodge structure on $S$.
  Here we notice that 
the Griffiths transversality of $H^{(j)}$ over $S^*$ clearly extends over $S$.
  Also, we remark that 
the canonical extension of Deligne (\cite{D70}) is nothing but 
the underlying $\cO_S$-module $\tau_*(H_{\bZ}^{(j)} \otimes \cO_S^{\log})$ of 
the extended $H^{(j)}$. 

For $s\in S$, we have the action $\rho_j:\pi_1(s^{\log}) \to \Aut(H_0)$ which is determined modulo $\G$.

\medskip

{\bf 5.3.}
Let $v$ be a point of $S_{\val}$ lying over $s$, and let $V_v$ be the valuative submonoid of $(M_S^{\gp}/\cO^\times_S)_s$ corresponding to $v$ containing $(M_S/\cO_S^\times)_s$. 
Let $\pi(v)$ be the $\bZ$-dual of $(M_S^{\gp}/\cO_S^\times)_s/V_v^\times$. 
If we log blow up $S$ sufficiently around $s$ and obtain $S'$ then, for the image $s'$ of $v$ in $S'$, there is a canonical isomorphism $\pi_1((s')^{\log})\simeq \pi(v)$ which is compatible with the isomorphism between $\pi_1(s^{\log})$ and the $\bZ$-dual of $(M_S^{\gp}/\cO_S^\times)_s$.

Note that both $\pi_1(s^{\log})$ and $\pi(v)$ are free abelian groups of finite rank. 

\bigskip

\head
\S6. Proof of Theorem 0.1
\endhead

{\bf 6.1.} 
Let the situation be as in 5.1.

We will prove (A) in 5.1 by induction on the log rank of $s$. 

By 5.2 and 4.6, 
we see that the map $f_j : S^*\to \G \bs D$ extends uniquely to a continuous map $S_{\val}^{\log}\to \G \bs D_{\val}^{\sharp}$ which we still denote by $f_j$. 
Let $Z=\{p\in S_{\val}^{\log}\;|\; f_1(p)=\dots=f_n(p)\}$. 
Let $p\in Z$ and assume $p\in s_{\val}^{\log}\subset S_{\val}^{\log}$ (that is, $p$ lies over $s\in S$). 
Since  $D_{\val}^{\sharp}\to \G \bs D_{\val}^{\sharp}$ is a local homeomorphism (2.3 (ii)), there is an open neighborhood $U_p$ of $p$ in $S^{\log}_{\val}$ such that, for $1\leq j\leq n$, $f_j : U_p\to \G \bs D_{\val}^{\sharp}$ lifts to a continuous map $\tilde f_j: U_p\to D_{\val}^{\sharp}$ and such that $Z_p:=\{p'\in U_p\;|\;\tilde f_1(p')=\dots =\tilde f_n(p')\} = U_p\cap Z$. 
We may assume that there is an open set $U'_p$ of $S^{\log}$ such that the image of $U_p$ in $S^{\log}$ is contained in $U'_p$ and the $H^{(j)}_\bZ$ are constant sheaves on $U'_p$. 

Let $Z_p(s)= Z_p\cap s_{\val}^{\log}$.
Since $s_{\val}^{\log}$ is compact, there is a finite subset $B$ of $s_{\val}^{\log}$ such that $\bigcup_{p\in B}\; Z_p(s)= Z\cap s_{\val}^{\log}$. 

We have either the following Case 1 in 6.2 or Case 2 in 6.3. 
For $p'\in S_{\val}^{\log}$ lying over $s$, let $\pi(p'):=\pi(v)$ (5.3), 
where $v$ is the image of $p'$ in $S_{\val}$. 

\medskip

{\bf 6.2.} 
{\it Case} 1. 
For some $p\in B$, the images of $\bQ\otimes \pi(p') \to \bQ\otimes \pi_1(s^{\log})$ for $p'\in Z_p(s)$ generate $\bQ\otimes \pi_1(s^{\log})$. 

\medskip

In this case, all $\rho_j$ ($1\leq j\leq n$) coincide. 
This is because, for any $p'\in Z_p(s)$, $\tilde f_1(p')=\dots =\tilde f_n(p')$ by definition, hence in particular the restrictions of $\rho_j$ to $\pi(p')$ for $1\leq j\leq n$ coincide.

In such a case, we can apply \cite{KNU.p} 7.4.1 (ii) to $f_j$ simultaneously. 
More precisely, 
let $\Gamma'$ be the image of $\pi_1(s^{\log})$ in $\Aut(H_0)$ under the common representations $\rho_j$ $(1\le j\le n)$.
Then, for every $j$, the period map $f_j: S^*\to \G\bs D$ lifts to $S^*\to \G' \bs D$. 

Let $C$ be the union over $1 \le j \le n$ of the sets of local monodromy cones of $H^{(j)}$ in $\fg_\bQ$ (cf.\ \cite{KU09} 2.5.11).
Using the argument in the proof of \cite{KU09} 4.3.8, we can construct a fan $\Sigma$ in $\fg_\bQ$ which satisfies the following conditions (1)--(3) (cf.\ ibid.\ 4.3.6).

\smallskip

\noindent
$(1)$ $\tCu_{\sig\in \Sig}\, \sig=\tCu_{\sig\in C} \,\sig$.
\medskip

\noindent
$(2)$ $\Sig$ is compatible with $\G'$.
\smallskip

\noindent
$(3)$ For any $\sig\in C$, $\sig=\bigcup_{j=1}^m\, \tau_j$ for some
$m \geq 1$ and for some $\tau_j\in \Sig$.
\smallskip

Then, similarly to ibid.\ 4.3.6, $\Sigma$ is strongly compatible with $\Gamma'$ and the liftings $S^*\to \G' \bs D$ of the period maps $f_j$ extend to $S(\Sigma)\to \G'\bs D_{\Sigma}$, where $S(\Sigma)$ is the log modification of $S$ corresponding to $\Sigma$ (ibid.\ 3.6).
From this and the properness of $S(\Sigma)\to S$, we have (A). 

\medskip

{\bf 6.3.} 
{\it Case} 2. 
For any $p\in B$, the images of $\bQ\otimes \pi(p') \to \bQ\otimes \pi_1(s^{\log})$ for $p'\in Z_p(s)$ do not generate $\bQ\otimes \pi_1(s^{\log})$. 

\medskip

In this case, for each $p\in B$, take a non-trivial $\bQ$-linear map
$h_p:\bQ \otimes \pi_1(s^{\log})\to \bQ$ which kills the images $\bQ\otimes \pi(p')$ for all $p'\in Z_p(s)$. 
For each $p\in B$, take the division of $\pi_1^+(s^{\log})$ by $A_p^{0}:=\{x\in \pi_1^+(s^{\log})\;|\;h_p(x)= 0\}$, 
$A_p^{+}:=\{x\in \pi_1^+(s^{\log})\;|\;h_p(x)\geq  0\}$, 
$A_p^{-}:=\{x\in \pi_1^+(s^{\log})\;|\;h_p(x)\leq 0\}$. 
Then we have clearly
\medskip

(P) For $p'\in Z_p(s)$, the image of $\bQ\otimes \pi(p')$ 
 in $\bQ\otimes \pi_1(s^{\log})$ is contained in $\bQ\otimes A_p^{0}$. 

\medskip

Define the subdivision of $\pi_1^+(s^{\log})$ consisting of all cones of the form $R\cap \bigcap_{p\in B} A_p^{c(p)}$, where $R$ is a face of $\pi_1^+(s^{\log})$ and $c(p)\in \{0, +, -\}$. 
Let $S'$ be the log modification of an open neighborhood of $s$ in $S$ corresponding to this subdivision. 
Let $S''$ be a log modification of $S'$ which is smooth. 

By (P), any point of the closure of $Z^*$ in $S''$ lying over $s$ is of log rank strictly smaller than the log rank of $s$. 

By the induction on the log rank, we can complete the proof of (A) and hence 
the proof of Theorem 0.1.

\enddocument